\documentclass[10pt,reqno]{amsart}
\usepackage{color, amsmath,amssymb, amsfonts, amstext,amsthm, latexsym}

\usepackage{verbatim}
\usepackage[utf8]{inputenc}

\newtheorem{Th}{Theorem}[section]
\newtheorem{theorem}[Th]{Theorem}
\newtheorem{assumption}{Assumption}
\newtheorem{definition}[Th]{Definition}
\newtheorem{corollary}[Th]{Corollary}
\newtheorem{proposition}[Th]{Proposition}
\newtheorem{lemma}[Th]{Lemma}

\let \ssection=\section
\renewcommand{\section}{\setcounter{equation}{0}\ssection}

\newcommand{\finX}[1]{\textbf{x}_{#1}}
\newcommand{\R}{\mathbb R}
\newcommand{\E}{\mathbb E}
\newcommand{\PP}{\mathbb P}
\newcommand{\N}{\mathcal N}
\newcommand{\LL}{\mathcal L}
\newcommand{\F}{\mathcal F}

\newcommand{\G}{C}
\newcommand{\Ksubdivision}{\kappa}

\newcommand{\dBr}{d\overleftarrow{B}_r}
\newcommand{\dB}[1]{d\overleftarrow{B}_{#1}}
\newcommand{\rv}{\right|}
\newcommand{\lv}{\left|}
\newcommand{\lp}{\left (}
\newcommand{\rp}{\right )}

\newcommand{\Fchapeau}{\widehat{\mathcal F}}

\newcommand{\iti}[2]{\int_{t_{#1}}^{t_{#2}}}
\newcommand{\schemeZhangY}[1]{Y^{\pi}_{t_{#1}}}
\newcommand{\schemeZhangZ}[1]{Z^{\pi}_{t_{#1}}}
\newcommand{\schemeZhangZcontinuous}[1]{Z^{\pi}_{#1}}
\newcommand{\schemeZhangZpi}[1]{Z^{\pi,1}_{t_{#1}}}
\newcommand{\schemeZhangX}[1]{X^{\pi}_{t_{#1}}}
\newcommand{\schemeZhangTheta}[1]{\Theta^{\pi,1}_{t_{#1}}}
\newcommand{\vectorX}[1]{\textbf{x}^{#1}}

\newcommand{\zpi}[1]{Z^{\pi}_{#1}}
\newcommand{\e}{\epsilon}
\newcommand{\tinitial}{0}
\newcommand{\Lf}{L_f}
\newcommand{\Lg}{L_g}
\newcommand{\lc}{\left [}
\newcommand{\rc}{\right ]}

\newcommand{\ssXt}[1][t]{\sigma\lp X_{#1}\rp}
\newcommand{\nnXt}[1][t]{\lp\nabla X_{#1}\rp^{-1}}
\newcommand{\ag}{\alpha}

\newcommand{\undeux}{\frac{1}{2}}

\newcommand{\EspaceY}{\R}
\newcommand{\EspaceZ}{\R}
\newcommand{\normeY}[1]{\lv #1\rv}
\newcommand{\normeZ}[1]{\lv #1\rv}
\newcommand{\normeYnabla}[1]{\lv #1\rv}

\newcommand{\normeX}[1]{\lv #1\rv}
\newcommand{\normeXnabla}[1]{\lv #1\rv}

\newcommand{\normeSIGma}[1]{\lv #1\rv}

\newcommand{\normeMatricekl}[1]{\lv #1\rv_{}}

\newcommand{\Z}{Z}
\newcommand{\W}{W}

\newcommand{\sumDelta}{}
\newcommand{\dW}[1][r]{dW^{}_{#1}}

\newcommand{\sumLambda}{}
\newcommand{\gl}{g}

\begin{document}
\title[Discretization of BDSDEs]
{On the discretization of backward doubly stochastic differential equations}
\author{Omar Aboura}
\address{ 
  SAMOS, Centre d'\'Economie de la Sorbonne,
 Universit\'e Paris 1 Panth\'eon Sorbonne,
90 Rue de Tolbiac, 75634 Paris Cedex France }
 \email{omar.aboura@malix.univ-paris1.fr}

\begin{abstract}
In this paper, we are dealing with the approximation of the process $\lp X_t,Y_t,Z_t\rp$ 
 solution to the backward doubly stochastic differential equation (BDSDE)
\begin{align*}
X_s&=x+\int_{\tinitial}^sb\lp X_r\rp dr+\int_{\tinitial}^s\sigma(X_r)\dW,\\
Y_{s}&=\phi\lp X_T\rp +\int_{s}^Tf\lp r,X_r,Y_r,Z_r\rp dr+\int_{s}^T g\lp r,X_r,Y_r,Z_r\rp d\overleftarrow{B}_r-\int_{s}^T Z_r\dW.
\end{align*}
After proving the $L^2$-regularity of $Z$, we use the Euler scheme to discretize $X$ 
and the  Zhang approach in order to give a discretization scheme  of the process $(Y,Z)$.
\end{abstract}
\maketitle

\section{Introduction}

Since the pioneering work of E.~Pardoux and S.~Peng  \cite{PPquasilinear}, 
 backward stochastic differential equations (BSDEs) have been intensively studied during the two last decades. 
Indeed, this notion has been a very useful tool to study problems in many areas, such as
mathematical finance, stochastic control, partial differential equations; see e.g.  \cite{MaYong}
where many applications are described.
Discretization schemes for BSDEs have been introduced and studied by several authors. 
The first papers on this  topic are that of V.Bally \cite{BallyDiscretization}  and D.Chevance \cite{ChevanceThesis}.
In his thesis, Zhang made an interesting contribution which was the starting point of intense study among
which the works of   B.~Bouchard and N.Touzi \cite{BouchardTouzi}, E.Gobet, J.P. Lemor and
X. Warin\cite{GobetRegression},...
The notion of BSDE has been generalized by E.~Pardoux and S.~Peng  \cite{PardouxPengSPDE} to that of 
 Backward Doubly Stochastic Differential Equation (BDSDE) as follows. Let  $(\Omega,\F,\PP)$ be
a  probability space,  $T$ denote some fixed terminal time which will be used throughout the paper, 
$\lp W_t\rp_{\tinitial\leq t\leq T}$ and $\lp B_t\rp_{\tinitial\leq t\leq T}$ be two independent standard
 Brownian motions defined on  $(\Omega,\F,\PP)$ and  with values in $\R^d$, and  $\R$ respectively. 
 On this space  we will deal with two families of $\sigma$-algebras:
\begin{equation}\label{defF}
 \F_t:=\F^W_{\tinitial,t}\vee\F_{t,T}^B\vee\N,\quad
\Fchapeau_t:=\F^W_{\tinitial,t}\vee\F_{\tinitial,T}^B\vee\N, \quad 
{\mathcal H}_t=\F^W_{\tinitial,T}\vee\F_{t,T}^B\vee\N,
\end{equation}
where $\F_{t,T}^B:=\sigma\lp B_r-B_t; t\leq r\leq T\rp$, 
$\F^W_{\tinitial,t}:=\sigma\lp W_r-W_{\tinitial};\tinitial\leq r\leq t \rp$ and
$\N$ denotes the class of $\PP$ null sets. We remark that  $(\Fchapeau_t)$ is a filtration, 
$({\mathcal H}_t)$ is a decreasing family of $\sigma$-albegras,  while 
$(\F_t)$ is neither  increasing nor decreasing. 
Given an initial condition  $x\in \R^d$,  let $(X_t)$ be the  $d$-dimensional diffusion process
 defined by
\begin{equation}\label{defX}
X_t = x+\int_{\tinitial}^tb\lp X_r\rp dr+\int_{\tinitial}^s\sigma\lp X_r\rp \dW .
\end{equation}
Let  $\xi\in L^2(\Omega)$ be an $\R^d$-valued, $\F_T$-measurable random variable, 
$f$ and $g$ be  regular enough coefficients; 
consider the BDSDE 
defined as follows: 
\begin{align} 
 Y_s\; = \; \xi & 
+\int_s^Tf\lp r,X_r,Y_r, Z_r\rp dr
\nonumber \\
&\quad 
+\sumLambda\int_s^T\gl\lp r,X_r,Y_r,Z_r\rp \dBr -\int_s^TZ_rdW_r.\label{BDSDE}
\end{align}
In this equation,  $dW$ is the forward integral and $\dB{}$ is the backward integral
 (we send the reader to \cite{NualartPardoux anticipating} for more  details on backward integration).
A solution to \eqref{BDSDE} is a  pair of  real-valued process $(Y_t,Z_t)$,
 such that $X_t$ and $Y_t$ are $(\F_t)$ for every $t\in [0,T]$, such that \eqref{BDSDE} is satisfied  and
\begin{equation} \label{bndYZ}
 \E\Big(\sup_{0\leq s\leq T} |Y_s|^2\Big) +  \E \int_0^T |Z_s|^2 ds<+\infty.
\end{equation}
In \cite{PardouxPengSPDE} Pardoux and Peng have proved that under some 
 Lipschitz property on $f$ and $g$ which will be stated more precisely in
 section \ref{s2}, 
\eqref{BDSDE} has a unique solution  $(Y,Z)$.

The aim of this paper is  to study the discretization of a Backward Doubly Stochastic Differential Equation
For the sake of simplicity, as in Zhang's paper \cite{Zhang},  we  assume  that 
$Y$ and $Z$ are real-valued processes. The extension to higher dimension is cumbersome 
and without theoretical problems. 
This discretization scheme of  $(Y,Z)$ is motivated by the  link between (\ref{BDSDE})
 and the following backward stochastic partial differential equation when $\xi=\phi(X_T)$ for a regular function $\phi$:
\begin{eqnarray}
u(t,x)&=&\phi(x)+\int_t^T\Big(\LL u(s,x)+f\lp s,x,u(s, x),\nabla u(s, x)\sigma (x)\rp \Big) ds\nonumber\\
&&
+\sumLambda\int_t^T\gl\lp s,x,u(s,x),\nabla u(s, x)\sigma (x)\rp \dB{s},\label{SPDE}
\end{eqnarray}
where $\LL$ is the differential operator defined by: 
$$
\LL u(t,x)=\undeux\sum_{i,j=1}^d\lp\sigma\sigma^*\rp_{i,j}(x)\frac{\partial^2}{\partial x_i\partial x_j}u(t,x)+\sum_{i=1}^db_i(x)\frac{\partial}{\partial x_i}u(t,x).
$$
The paper is organized as follows: first we prove the $L^2$-regularity of $Z$ in section \ref{s2}.
 This a crucial  step in order to  the scheme using Zhang's method, 
 which is done in section \ref{s3}. Finally,  a numerical scheme is described in the last section.
To ease notations, we  set $\Theta_r:=\lp X_r,Y_r,Z_r\rp$ for $r\in [0,T]$. As usual, we 
 denote by $C_p$  a constant which depends on some parameter $p$, and which can change
from on line to the next one. Finally, for some function $h(t,x,y,z)$ defined on 
$[0,T]\times \R^d\times \R\times \R$, we let $\partial_y h(t,x,y,z)$ (resp. $\partial_z h(t,x,y,z)$)
the partial derivatives of $h$ with respect to the real variable $y$ (resp.  $z$),
 while $\partial_x h(t,x,y,z)$
will denote the vector $(\partial_{x_i}h(t,x,y,z), i=1, \cdots, d)$.

\section{Regularity properties}\label{s2}
In this section we give some regularity properties of the process $X,Y$ and $Z$.

The following  assumptions which ensure existence and uniqueness of the solution will be in
 force throughout
  the paper.
For every integer $n\geq 1$, let $M^2([0,T],\R^n)$ denote the set of $\R^n$-valued jointly measurable processes
$(\varphi_t , t\in [0,T])$ such that $\varphi_t$ is $\F_t$-measurable for almost every $t$ and 
$\E\int_0^T |\varphi_t|^2 dt <+\infty$. 
\begin{assumption}[for the forward process $X$]\label{assumptionX}
The maps $b:\R^d\rightarrow\R^d$ and $\sigma: \R^d\rightarrow \R^{d\times d}$ 
are of class $\mathcal C^3_{b}$.
\end{assumption}
\begin{assumption}[for the backward process $(Y,Z)$]\label{assumptionYZ}
 Let 
$ f:\Omega\times [\tinitial,T]\times\R^d\times\EspaceY\times \EspaceZ  \rightarrow\EspaceY$, 
$\gl:\Omega\times [\tinitial,T]\times\R^d\times\EspaceY\times \EspaceZ  \rightarrow\EspaceY $
 be such that $f$ and $g$ are jointly measurable, for every 
$(x,y,z)\in \R^{d+2}$, $f(.,x,y,z)$ and $g(.,x,y,z)$ belong to $M^2([0,T],\R)$, and such that:
\begin{itemize}
\item[(i)] There exist some nonnegative constants  $\Lf,\Lg$ and  a constant
 $\ag\in[0,1)$ such that for every $\omega\in \Omega$,  $t,t'\in\lc  \tinitial,T\rc $,  $x,x'\in \R^d$, 
$y,y'\in\R$ and  $z,z'\in\R$
\begin{eqnarray*}
\normeY{ f\lp t,x,y,z \rp-f\lp t',x',y',z' \rp}^2&\leq &
\Lf \lp \lv t-t'\rv +\normeX{ x-x'}^2 +\normeY{ y-y'}^2+\normeZ{ z-z'}^2 \rp,\\
\normeY{\gl\lp t,x,y,z \rp-\gl\lp t',x',y',z' \rp}^2&\leq &\Lg \lp \lv t-t'\rv +\normeX{ x-x'}^2 +\normeY{ y-y'}^2 \rp
+\ag\normeZ{ z-z'}^2,
\end{eqnarray*}
\item[(ii)] For all $s\in\lc  \tinitial,T\rc $ $f\lp s,.\rp$ and
 $g(s,.)$ are of class ${\mathcal C}^3$ with bounded partial derivatives up to order 3, uniformly in time.
\item[(iii)] 
For a function $h(t,x,y,z)$, set $h(t,0):=h(t,0,0,0)$. Then 
$$
\sup_{r\in[\tinitial,T]}\normeY{ f(r,0)}+\sup_{r\in[\tinitial,T]}\normeMatricekl{ g(r,0)}<\infty.
$$
\end{itemize}
\end{assumption}
\begin{assumption}
Suppose that $\xi:= \phi(X_T)$ for some function $\phi:\R^d\to \R$ of class ${\mathcal C}^2_b$
and that 
 for every $\omega\in \Omega$, 
$$ \sup_{t,x,y,z}\lv {\partial_z} g\lp t,x,y,z\rp\rv < 1 .
$$
\end{assumption}

\subsection{Some classical properties of the forward process $X$}
We at first recall  without proof the following well known results  on diffusion processes.
Define the  
$\R^{d\times d}$-valued process $\lp\nabla X_t\rp_{\tinitial\leq t\leq T}$ by:

$$
\nabla X_t:=\Big(\frac{\partial}{\partial x_j}X_t^i , i,j=1, \cdots, d\Big).   
$$

Then $\nabla X_t$ is an invertible $d\times d$  matrix, 
 solution to a linear stochastic differential equation with  coefficients depending on $X_t$.
Furthermore, the assumptions on the coefficients $\sigma$ and $b$ yield the  following classical result:
\begin{proposition}\label{propositionForwardX}
(i) For all $p\geq1$, there exist a constant $C_p>0$ such that for all $t,s\in\lc \tinitial,T\rc $:
 $$
\E \normeX{X_t -X_s }^{2p} + \E\normeXnabla{\nnXt -\nnXt[s] }^{2p}
\leq C_p|t-s|^p. 
$$
(ii) For all $p\in [1,+\infty[$, there exist a constant $C_p>0$ such that 
$$
\E\Big( \sup_{t\in[\tinitial,T]} |X_t|^{2p} +  \sup_{t\in[\tinitial,T]}\normeXnabla{\nnXt  }^p\Big)
\leq C_p. 
$$
\end{proposition}

\subsection{Time increments of $Y$ and $L^2$-regularity of $Z$}

The following lemma provides upper bounds for time increments of $Y$.
\begin{lemma}\label{lemmaYYZ}
 Set $\xi=\phi(X_T)$ for some  function $\phi:\R^d\to \R$ be of class ${\mathcal C}^1_b$.
 Then we have\\
(i) For all $p\geq 2$, there exist a constant $C_p>0$ depending on $T$
 such that for all $t,s\in\lc \tinitial,T\rc $
\begin{equation}\label{AccroisementY}
\E\normeY{Y_t-Y_s}^p\leq 
C_p|t-s|^{\frac{p}{2}} . 
\end{equation}
(ii) For all $p\geq 1$, there exist a constant $\G>0$ such that 
\begin{equation}\label{SupZ}
\sup_{\tinitial\leq r\leq T}\E \normeYnabla{ Z_r}^{2p}\leq C .
\end{equation}
\end{lemma}

\noindent Notice that the  inequality \eqref{AccroisementY}
 is  different from  equation (2.11) in \cite{Zhang}.

\begin{proof}  
  We at first prove (ii).
Let  $\lp\nabla Y_t\rp_{\tinitial\leq t\leq T}=
({\partial_x}Y_t)_{0\leq t\leq T}$  denote the real-valued process defined by differentiation of $Y$ 
 as function  of  the initial condition $x$ of the
diffusion process $(X_t)$.
We recall the following representation of $Z$ (see \cite{PardouxPengSPDE} Proposition 2.3):
\begin{equation}\label{reresentationZ}
Z_t=\nabla Y_t\nnXt\ssXt.
\end{equation}
where $\lp \nabla Y_t,\nabla Z_t\rp$ satisfies the linear BDSDE with the forward process $(X_t,\nabla X_t)$ and
the evolution equation:
\begin{align}\label{nablaY}
\nabla Y_t= & \phi'(X_T)\nabla X_T+\int_t^T\Big( f_x(r, \Theta_r)\nabla X_r+f_y(r, \Theta_r)\nabla Y_r
+f_z(r, \Theta_r)\nabla Z_r\Big) dr \nonumber \\
&+\sumLambda\int_t^T\Big(  \gl_x(r,\Theta_r)\nabla X_r+\gl_y(r,\Theta_r)\nabla Y_r
+\gl_z(r,\Theta_r)\nabla Z_r\Big)  \dBr  -\sumDelta\int_t^T\nabla \Z_r d\W_r .
\end{align}
By E.Pardoux and S.Peng \cite{PardouxPengSPDE} page 217, we deduce 
\begin{equation}\label{n=bndnablaY}
\E\Big( \sup_{\tinitial\leq t\leq T}\normeYnabla{ \nabla Y_t}^p\Big) <\infty.
\end{equation}
Then H\"older's inequality and Proposition \ref{propositionForwardX}
yield  
$$
\E\normeYnabla{ Z_t}^{2p}\leq \lp\E\normeYnabla{ \nabla Y_t}^{6p}\rp^{\frac{1}{3}}\lp\E\normeXnabla{ \lp \nabla X_t\rp^{-1}}^{6p}\rp^{\frac{1}{3}}\lp\E\normeSIGma{ \sigma \lp X_t\rp}^{6p}\rp^{\frac{1}{3}}.
$$
This concludes  the proof of (ii).

(i) Suppose that $s<t$, then using \eqref{BDSDE}, we deduce that
\begin{eqnarray*}
\normeY{Y_t-Y_s}^p
& \leq &
 C_p \normeY{\int_s^tf\lp r,\Theta_r\rp -f\lp r,X_r,Y_r,0\rp dr}^p+C_p \normeY{\int_s^tf\lp r,X_r,Y_r,0\rp dr}^p\\
&&+C_p\sumLambda\normeY{\int_s^t\gl(r,\Theta_r)\dBr}^p+C_p\sumDelta\normeY{\int_s^t\Z_rd\W_r}^p. 
\end{eqnarray*}
Recall that $\widehat{\F}_t$ and ${\mathcal H}_t$ have been defined in \eqref{defF}.
The process $(\int_0^t Z_r dW_r , 0\leq t\leq T)$ is a $(\widehat{\F}_t)$-martingale, while 
the process $(\int_t^T \gl(r,\Theta_r)\dBr ,0\leq t\leq T)$ is a backward martingale for 
$({\mathcal H}_t)$. Hence, the  Burkholder-Davies-Gundy and 
 H\"{o}lder inequalities yield  
\begin{eqnarray}
\E\normeY{ Y_t-Y_s}^p
& \leq&
 C_p|t-s|^{p-1}\E\int_s^t\normeY{f\lp r,X_r,Y_r,0\rp}^pdr
+ C_p \E\sumDelta\normeY{\int_s^t\lv \Z_r \rv dr}^p\nonumber\\
& &+C_p\sumLambda\E\lp\int_s^t\normeY{\gl(r,\Theta_r)}^2dr\rp^{\frac{p}{2}}
+C_p\sumDelta\E\lp\int_s^t\normeY{\Z_r}^2dr\rp^{\frac{p}{2}}.\label{YmoinsY1}
\end{eqnarray}
Assumption \ref{assumptionYZ} (i) and (ii), Proposition \ref{propositionForwardX}
 and \eqref{bndYZ} yield
\begin{align}
 \E\Big(\int_s^t  & \normeY{\gl(r,\Theta_r)}^2  dr\Big)^{\frac{p}{2}}
\leq
 C_p\E\lp\int_s^t\normeY{ \gl(r,\Theta_r)-\gl(r,0)}^2dr\rp^{\frac{p}{2}}
+C_p\E\lp\int_s^t\normeY{\gl(r,0)}^2dr\rp^{\frac{p}{2}}\nonumber\\
\leq&
C_p|t-s|^{\frac{p}{2}}  
+C_p\E\lp\int_s^t \big( \normeX{X_r}^2 + |Y_r|^2 \big) \, dr\rp^{\frac{p}{2}}
+C_p\sumDelta\E\lp\int_s^t\normeY{ \Z_r}^2dr\rp^{\frac{p}{2}}\nonumber\\
\leq&
C_p|t-s|^{\frac{p}{2}}+C_p\sumDelta\E\lp\int_s^t\normeY{ \Z_r}^2dr\rp^{\frac{p}{2}}\label{YmoinsY2}. 
\end{align}
Similarly,
\begin{align}
 \E\int_s^t  | f(r,&X_r,Y_r,0)|^pdr
\leq  C_p\E \int_s^t\normeY{ f( r,0)}^pdr  + 
 C_p\E\int_s^t\normeY{ f\lp r,X_r,Y_r,0\rp-f\lp r,0\rp}^pdr
\nonumber \\
\leq &
C_p\E \int_s^t \normeY{f\lp r,0\rp}^pdr+C_p \int_s^t \E\lp\normeX{X_r}^p+\normeY{Y_r}^p\rp dr 
\leq C_p\, |t-s|\label{YmoinsY3}  .
\end{align}
Hence, the inequalities (\ref{YmoinsY1})-(\ref{YmoinsY3}) 
imply 
\[ 
\E\normeY{Y_t-Y_s}^p\leq 
C_p|t-s|^{\frac{p}{2}}+C_p\sumDelta\E\lp\int_s^t\normeY{\Z_r}^2dr\rp^{\frac{p}{2}}.
\] 
Using H\"older's inequality and \eqref{SupZ} we conclude the proof of  \eqref{AccroisementY}.
\end{proof}
Since equation  \eqref{nablaY} proves that the pair $(\nabla Y, \nabla Z)$ is the solution of a 
BDSDE with forward process $(X, \nabla X)\in L^p$ for every $p\in [1,+\infty[$, we deduce from 
\eqref{AccroisementY}   that for every function $\phi:\R^d\to \R$ of class ${\mathcal C}^2_b$, 
we have for  $0\leq s<t\leq T$ and $p\in [1,+\infty[$: 
\begin{equation} \label{acroisnablaY}
\E| \nabla Y_t-\nabla Y_s|^p\leq 
C_p|t-s|^{\frac{p}{2}}  , 
\end{equation} 
 for some constant $C_p>0$.
We now  establish some control of time increments of the process  $Z$, 
 following the idea of   J.Zhang \cite{Zhang}. 

\begin{theorem}[$L^2$-regularity of $Z$]\label{L2regularity}  There exists  
a non negative constant $C$ such that for every subdivision $\pi=\{t_0=0<t_1\cdots <t_n=T\}$ with
mesh $|\pi|$, one has 
\begin{equation}\label{L2regu}
 \sum_{1\leq i\leq n}\E \int_{t_{i-1}}^{t_i}\lp\normeZ{ Z_t-Z_{t_{i-1}}} ^2+\normeZ{Z_t-Z_{t_i}}^2\rp dt \leq \G |\pi|.
\end{equation}

\end{theorem}
\proof
Using the representation of $Z$ as a product, we deduce  (\ref{reresentationZ}),
\begin{align*}
Z_t-Z_{t_i}=\nabla Y_t\lp \nabla X_t\rp^{-1}\sigma\lp X_t\rp-\nabla Y_{t_i}\lp \nabla X_{t_i}\rp^{-1}\sigma\lp X_{t_i}\rp.
\end{align*}
Then, 
\begin{align*}
\normeZ{Z_t-Z_{t_i}}^2\leq &
3\normeYnabla{\nabla Y_t-\nabla Y_{t_i}}^2\normeXnabla{\lp \nabla X_{t}\rp^{-1}}^2\normeXnabla{\sigma\lp X_{t}\rp}^2\\
&+3\normeYnabla{\nabla Y_{t_i}}^2
\normeXnabla{\lp \nabla X_t\rp^{-1}-\lp \nabla X_{t_i}\rp^{-1}}^2\normeXnabla{\sigma\lp X_{t}\rp}^2\\
&+3\normeYnabla{\nabla Y_{t_i}}^2\normeXnabla{\lp \nabla X_{t_i}\rp^{-1}}^2\normeXnabla{\sigma\lp X_t\rp-\sigma\lp X_{t_i}\rp}^2.
\end{align*}
To conclude the proof, we use H\"older's inequality,  Proposition \ref{propositionForwardX} and \eqref{acroisnablaY}.
\endproof
Theorem \ref{L2regularity} immediatly yields the following
\begin{corollary}\label{corZ}
$$
\sum_{1\leq i\leq n-1}\E\int_{t_{i-1}}^{t_{i+1}}|Z_r-Z_{t_i}|^2dr\leq C|\pi|.
$$
\end{corollary}

\section{The  discretization of $(X,Y,Z)$}\label{s3}
\subsection{Discretization of the process $X$: The Euler scheme}
We briefly recall the Euler scheme and send the reader to \cite{KlodenPlaten} for more details.
Let  $\pi:=\{t_0=\tinitial<t_1<...<t_n=T\}$ be a subdivision of $[0,T]$. 
 We define the process $X^{\pi}_{t}$, called the Euler scheme, by
$$
X^{\pi}_{t}=X^{\pi}_{t_0}+\int_{t_0}^tb\lp X^{\pi}_{s_{\pi}} \rp ds+\int_{t_0}^t\sigma(X^{\pi}_{s_{\pi}}) dW_s,
$$
where $s_{\pi}:=\max\{ t_i\leq s\}$.
The following result is well known: 
\begin{proposition}\label{propositionX} There exists a constant $C>0$ such that for every subdivision $\pi$,
 $$\max_{i}\E\normeX{ X_{t_i}-X^{\pi}_{t_i}}^2\leq \G|\pi|,\qquad
\E\iti{i-1}{i} \normeX{ X_r-X^{\pi}_{t_i}}^2 dr\leq C|\pi|^2.
$$
\end{proposition}
\subsection{Discretization of the process $(Y,Z)$: The step process}
In this section, we construct an approximation of $(Y,Z)$  using Zhang's approach.\\
Let $\pi:t_0=\tinitial<...<t_n=T$ be 
any subdivision 
on $\lc \tinitial,T\rc $. 
 Set ${\mathcal G}_t={\mathcal G}^i_t$ for  $t_{i-1}\leq t <  t_i$,
where we let  
$$\mathcal G^i_t:=\sigma\lp W_r-W_{\tinitial};\tinitial\leq r\leq t\rp
\vee\sigma\lp B_r-B_{t_{i-1}};t_{i-1}\leq r\leq T  \rp\, ,\quad t_{i-1}\leq t \leq  t_i,  $$
 and define the $(\mathcal G_t)$-adapted process $\lp\schemeZhangY{},\schemeZhangZ{}\rp_{0\leq t\leq T}$ 
recursively (in a backward manner), as follows:
Set $\schemeZhangY{n}=\phi(\schemeZhangX{n}),\;\schemeZhangZpi{n}=0$; 
for  $i=n-1,...,0$, let 
$$
\schemeZhangZpi{i}:=\frac{1}{\Delta
  t_{i+1}}\E\lp\left.\iti{i}{i+1}\schemeZhangZcontinuous{r}dr\right|\F_{t_i}\rp,
$$ 
 
and for  $i=n,...,1$, let 
$$
\Delta t_i=t_i-t_{i-1}, \Delta B_{t_i}=B_{t_i}-B_{t_{i-1}},
\Theta_{t_i}^{\pi,1}:=
\lp \schemeZhangX{i},\schemeZhangY{i},\schemeZhangZpi{i}\rp,
$$
\begin{equation}\label{SchemeZhang}
\schemeZhangY{}=\schemeZhangY{i}+ f\lp t_i,\Theta_{t_i}^{\pi,1}\rp \Delta t_i
+ g\lp t_i,\Theta_{t_i}^{\pi,1}\rp \Delta B_{t_i}-\iti{}{i} Z_r^{\pi}dW_r ,\quad \forall  t\in[t_{i-1},t_i).
\end{equation}
 Note that the equation (\ref{SchemeZhang}) is not a BDSDE in the sense of \cite{PardouxPengSPDE}; 
however,  we have the following:
\begin{proposition}\label{ZhangExistence}
For every $i=1,...,n$, there exists a process $(\schemeZhangY{},\schemeZhangZ{})_{t\in[t_{i-1},t_i)}$
 adapted to the filtration $({\mathcal G}_t , t_{i-1}\leq t<t_i)$, 
 such that (\ref{SchemeZhang}) holds. Furthermore,
$\schemeZhangY{i}\in\F_{t_i}$.
\end{proposition}
\proof
The proof is similar to that 
in \cite{PardouxPengSPDE} page 212 and relies on the martingale representation theorem.
Fix an integer $i>0$ and suppose that the processes $(Y^{\pi}_t)$ and $(Z^{\pi}_t)$ have been defined for
$t\geq t_i$,  $({\mathcal G}_t)$-adapted,  and that $Y^\pi_{t_k}$ is $\F_{t_k}$-measurable for
$k=i, \cdots, n$.   We denote by 
$\lp M^i_t\rp _{t\in [t_{i-1},t_i] }$ the process defined by
$$
M^i_t:=\E\lp \schemeZhangY{i}+f\lp\left. t_i,\schemeZhangTheta{i}\rp\Delta t_i+g\lp t_i,
 \schemeZhangTheta{i}\rp\Delta B_{ t_i} \rv\mathcal G^i_t\rp, \quad t_{i-1}\leq t\leq t_i.
$$
By the martingale representation theorem, there exists  a $\lp\mathcal
G^i_t,t_{i-1}\leq t\leq t_i\rp$-adapted and square 
integrable process $\lp N_t^i, t_{i-1}\leq t\leq t_i\rp $ such that for $t_{i-1}\leq t\leq t_i$, 
$M^i_t=M^i_{t_{i-1}}+\iti{i-1}{}N^i_sdW_s$.
Therefore, 
$M^i_t=M^i_{t_i}-\iti{}{i}N_s^i dW_s$.
Clearly, 
$\mathcal G^i_{t_i} $ contains $\F_{t_i}$, $X^\pi_{t_i}$ is $\F_{t_i}^W\subset \F_{t_i}$ measurable
  and $\Theta_{t_i}^{\pi,1}$ is $\F_{t_i}$-measurable; hence  
$$
M^i_{t_i}=\schemeZhangY{i}+ f\lp t_i,\Theta_{t_i}^{\pi,1}\rp \Delta t_i+ g\lp t_i,\Theta_{t_i}^{\pi,1}\rp  \Delta B_{t_i}.
$$
Furthermore, note that ${\mathcal G}^i_{t_{i-1}} = \F_{t_{i-1}}$, so that $M^i_{t_{i-1}}$ is $\F_{t_i}$-measurable. 
This completes the proof by setting: $\schemeZhangY{}=M^i_t$, $\schemeZhangZ{}=N^i_t$ for $t_{i-1}\leq t<t_i$.
\endproof
Before stating the main theorem  of  this section, we introduce the following
\begin{definition}
 Let $\Ksubdivision\geq1$ be a constant. The subdivision $\pi$ is said to be
 $\Ksubdivision$-uniform if $\Ksubdivision\Delta t_i\geq |\pi|$ for every $i\in\{1,...,n\}$ .
\end{definition}
The main example of a $\Ksubdivision$-uniform subdivision is a uniform subdivision
 (i.e. for all $i$, $\Delta t_i=|\pi|$) where $\Ksubdivision=1$.
The following lemma gives an upper estimate of $Z_{t_i}-Z^{\pi,1}_{t_i} $.
\begin{lemma}\label{lem}For any $i=0,...,n-1$,  any $\kappa$-uniform subdivision $\pi$ and  $\beta>0$ we have:
 $$
\Delta t_i \E\normeZ{ Z_{t_i}-Z^{\pi,1}_{t_i}}^2
\leq
 \Ksubdivision\lp1+\beta\rp\int_{t_i}^{t_{i+1}}\normeZ{\zpi{r}-Z_r}^2dr
+
\Ksubdivision\lp1+\beta^{-1}\rp\int_{t_i}^{t_{i+1}}\normeZ{ Z_r-Z_{t_i}}^2dr.
$$
\end{lemma}
\proof For any $i=0,...,n-1$, $Z_{t_i}$ is $\F_{t_i}$-measurable, and  $\Delta t_i\leq |\pi|\leq\kappa \Delta t_{i+1}$; thus
\begin{align*}
\Delta t_i&  \E\normeZ{ Z_{t_i}-Z^{\pi,1}_{t_i}}^2=\Delta t_i \E\normeZ{ Z_{t_i}-\frac{1}{\Delta t_{i+1}}
\E\lp \left.\int_{t_i}^{t_{i+1}}\zpi{r}dr \rv\F_{t_i}\rp }^2\\
&
=\frac{\Delta t_i}{\lp\Delta t_{i+1}\rp^2}\E\normeZ{\E\lp 
\left.\int_{t_i}^{t_{i+1}}\lp Z_{t_i}-\zpi{r}\rp dr \rv\F_{t_i}\rp }^2\\
&
\leq\frac{\Ksubdivision}{\Delta t_{i+1}}\E\normeZ{\int_{t_i}^{t_{i+1}}\lp Z_{t_i}-\zpi{r}\rp dr}^2
\leq \Ksubdivision\E\int_{t_i}^{t_{i+1}}\normeZ{ Z_{t_i}-\zpi{r}}^2dr.
\end{align*}
where the last step is deduced from Schwarz's inequality. Using the usual estimate $|Z_{t_i}-Z_r^\pi|^2
\leq (1+\beta) |Z_r^\pi - Z_r|^2 + (1+ \beta^{-1}) |Z_r-Z_{t_i}|^2$, we conclude the proof.
\endproof
The following theorem is the main result of this section. It proves that as
 $|\pi|\rightarrow0$, $\lp Y^{\pi},Z^{\pi}\rp$ converges to $\lp Y,Z\rp$.
\begin{theorem}\label{th}
 Let $\pi$ be a $\Ksubdivision$-uniform subdivision with  sufficiently small mesh $|\pi|$, 
   $\ag<\frac{1}{\Ksubdivision}$,  let $\phi\in {\mathcal C}^2$ and $\xi=\phi(X_T)$. Then we have
\begin{equation}\label{ZhangError}
 \max_{0\leq i\leq n}\E\normeY{ Y_{t_i}-\schemeZhangY{i}}^2
+\E\int_{\tinitial}^T\normeZ{ Z_r-\schemeZhangZcontinuous{r}}^2dr
\leq
 C|\pi|. 
\end{equation}
\end{theorem}
\proof
Set $I_n=\E |\phi(X_T)-\phi(X^\pi_T)|^2$ and for  $i=1,...,n$, let 
$$
I_{i-1}:=\E \normeY{ Y_{t_{i-1}}-\schemeZhangY{i-1}}^2+\E\iti{i-1}{i} \normeZ{ Z_r-Z_r^{\pi}}^2dr.
$$
Using (\ref{BDSDE}) with $\xi=\phi(X_T)$ and (\ref{SchemeZhang}),  we deduce
\begin{align}\label{difY}
Y_{t_{i-1}}-Y_{t_{i-1}}^{\pi}+\iti{i-1}{i} \lp Z_r-Z_r^{\pi}\rp dW_r= &
Y_{t_{i}}-Y_{t_{i}}^{\pi}+\iti{i-1}{i} \lp f\lp r,\Theta_r\rp-f\lp t_i,\schemeZhangTheta{i}\rp\rp dr 
\nonumber \\
& +\iti{i-1}{i} \lp g\lp r,\Theta_r\rp-g\lp t_i,\schemeZhangTheta{i}\rp\rp d\overleftarrow{B}_r .
\end{align}
By construction, $Y_{t_{i-1}}-\schemeZhangY{i-1}$ is $\F_{t_{i-1}}=\mathcal G^i_{t_{i-1}}$
 measurable while for $r\in [t_{i-1},t_i)$,  $Z_r-Z^{\pi}_{r}$ is $(\mathcal G_r)$-adapted. 
 Hence, $Y_{t_{i-1}}-\schemeZhangY{i-1}$ is orthogonal to $\iti{i-1}{i}\lp Z_r-Z^{\pi}_{r}\rp dW_r$.
Therefore,
\begin{align*}
 I_{i-1}=&\E\normeY{ Y_{t_{i-1}}-Y_{t_{i-1}}^{\pi}+\iti{i-1}{i} \lp Z_r-Z_r^{\pi}\rp dW_r}^2 .
\end{align*}
Since   $g(r,\Theta_r)$ (resp. $g(t_i,\Theta_{t_i}^{\pi,1})$) is $\F_r$ ( resp. $\F_{t_i}$)-measurable,
the random variables  $Y_{t_{i}}-Y_{t_{i}}^{\pi}$ and
 $\iti{i-1}{i}\lp  g(r,X_r,Y_r)-g\lp t_i,\schemeZhangX{i},\schemeZhangY{i}\rp\rp \dBr$ 
are orthogonal. Hence for every $\e>0$, using assumption \ref{assumptionYZ},
 the $L^2$-isometry  of backward stochastic integrals,  Schwarz's inequality
 and \eqref{difY}, we deduce
\newcommand{\yy}{\normeY{ Y_{t_{i}}-Y_{t_{i}}^{\pi}}^2}
\begin{align*}
 I_{i-1}& \leq \lp 1+\frac{\Delta t_i}{\e}\rp\E\yy+\lp 1
+
2\frac{\e}{\Delta t_i} \rp\E\normeY{\iti{i-1}{i}\lp f\lp r,\Theta_r\rp-f\lp t_i,\schemeZhangTheta{i}\rp\rp dr}^2\\
&
+
\lp 1+\frac{\Delta t_i}{\e}\rp\E\normeY{\iti{i-1}{i}\lp g\lp r,\Theta_r\rp-g\lp t_i,\schemeZhangTheta{i}\rp\rp 
d\overleftarrow{B}_r}^2\\
\leq &
\lp 1+{\Delta t_i}{\e}^{-1} \rp\E\yy+\lp \Delta t_i+2\e \rp\E\iti{i-1}{i}\normeY{ f\lp r,\Theta_r\rp-f\lp t_i,
\schemeZhangTheta{i}\rp}^2 dr\\
&
+\lp 1+{\Delta t_i}{\e}^{-1} \rp\E\iti{i-1}{i}\normeMatricekl{ g\lp r,\Theta_r\rp
-g\lp t_i,\schemeZhangTheta{i}\rp}^2 dr\\
\leq &
\Big[ 1+{\Delta t_i}{\e}^{-1} +2\Lf\lp \Delta t_i^2+2\e\Delta t_i\rp
+2\Lg\lp \Delta t_i +{\Delta t_i^2}{\e}^{-1} \rp \Big] \, \E\yy\\
&
+\Big[ \Lf\lp \Delta t_i+2\e \rp+\Lg\lp 1+{\Delta t_i}{\e}^{-1} \rp\Big]
 \E\iti{i-1}{i}\!\! \lp|\pi|+\normeX{X_r-X^{\pi}_{t_i}}^2+2\normeY{ Y_r-Y_{t_i}}^2 \rp dr\\
&
+\Big[ \Lf\lp \Delta t_i+2\e \rp+\ag\Big( 1+{\Delta t_i}{\e}^{-1} \Big) \Big]
\E\iti{i-1}{i}\normeZ{ Z_r-\schemeZhangZpi{i}}^2 dr .
\end{align*}
For $|\pi|\leq 1$, $\Delta t_i^2\leq \Delta t_i$; using Proposition \ref{AccroisementY} with $p=2$ and Proposition
  \ref{propositionX},  we deduce
$$
\E\iti{i-1}{i}\lp |\pi|+\normeX{X_r-X^{\pi}_{t_i}}^2+2\normeY{ Y_r-Y_{t_i}}^2 \rp dr\leq C|\pi|^2 ,
$$
for some constant $C>0$. Hence for any $\gamma>0$
\begin{align*}
 I_{i-1}\leq &\Big[  1+ \Big( \e^{-1} +2\Lf(1+2\e) +2\Lg \big(1+\e^{-1}\big)  \Big) \Delta t_i \Big]\, \E\yy\\
&+C\Big[ \Lf\lp \Delta t_i+\e \rp+\Lg\lp 1+{\Delta t_i}{\e}^{-1} \rp\Big]\, |\pi|^2\\
&
+\big( 1+{\gamma}^{-1}\big)  \big[  \Lf\lp \Delta t_i+2\e \rp+\ag\lp 1+{\Delta t_i}{\e}^{-1} \rp\big]\, 
\E\iti{i-1}{i} \normeZ{ Z_r-Z_{t_i}}^2 dr\\
&
+\lp1+\gamma\rp\big[ \Lf\lp \Delta t_i+2\e \rp+\ag\lp 1+{\Delta t_i}{\e}^{-1} \rp\big] 
\Delta t_i\E\normeZ{ Z_{t_i}- 
\schemeZhangZpi{i}}^2.
\end{align*}
\newcommand{\ccce}{C_{\e}}
Lemma \ref{lem} yields for some positive constants $\ccce$, $C_{\e,\gamma}$ and $C_{\e, \gamma,\beta}$, we have:
\begin{align} \label{ETOILE}
&  I_{i-1}\leq \lp 1+ \ccce\Delta t_i \rp\E\yy+\ccce \, |\pi|^2
+C_{\e,\gamma}\; \E\iti{i-1}{i}\normeZ{ Z_r-Z_{t_i}}^2 dr\nonumber\\
&\quad 
+\Ksubdivision\lp 1+\gamma\rp  \lp 1+\beta\rp  \Big[ \Lf\lp \Delta t_i+2\e \rp
+\ag\lp 1+{\Delta t_i}{\e}^{-1} \rp \Big] \E 
\iti{i}{i+1} \normeZ{ Z^{\pi}_r-Z_r}^2dr\nonumber\\
&\quad 
+\Ksubdivision\lp 1+{\gamma} \rp  \lp 1+{\beta}^{-1}\rp  \Big[ \Lf\lp \Delta t_i+2\e \rp
+\ag\lp 1 +{\Delta t_i}{\e}^{-1} \rp \Big] \E \iti{i}{i+1} \normeZ{ Z_r-Z_{t_i}}^2dr
\nonumber\\
&\leq 
\lp 1+ \ccce\, \Delta t_i \rp\E\yy+\ccce|\pi|^2
+C_{\e,\gamma,\beta}\; \E\iti{i-1}{i+1} \normeZ{Z_r-Z_{t_i}}^2 dr
\nonumber \\
&
\quad + \Ksubdivision\lp1+\gamma\rp\lp 1+\beta\rp 
\Big[ \Lf\lp \Delta t_i+2\e \rp+\ag\lp 1+{\Delta t_i}{\e}^{-1} \rp\Big] \; 
\iti{i}{i+1} \normeZ{ Z^{\pi}_r-Z_r}^2dr. 
\end{align}
Recall that $\ag<\frac{1}{\Ksubdivision}$ and let $0<\delta<1-\Ksubdivision\ag$. 
Then choose positive constants $\beta$ and $\gamma$ small enough to ensure 
$\Ksubdivision\lp1+\gamma\rp\lp1+\beta\rp\ag<
1-\frac{2\delta}{3}$. Finally, let $\e>0$ small enough to ensure that 
$2\Ksubdivision\lp1+\gamma\rp\lp1+\beta\rp\Lf\e<\frac{\delta}{6}$. Then (\ref{ETOILE}) implies 
the existence of $C>0$ such that
 for every $i=1,...,n-1$,  
\begin{equation}\label{cinqquatre}
 I_{i-1}+\frac{\delta}{3}\E \iti{i}{i+1}\!\!  \normeZ{ Z^{\pi}_r-Z_r}^2dr \leq \lp 1+ C \Delta t_i \rp I_i+C|\pi|^2
+C\E\iti{i-1}{i+1}\!\! \normeZ{ Z_r-Z_{t_i}}^2 dr.
\end{equation}
Using the discrete Gronwall lemma in \cite{Zhang}   (Lemma 5.4 page 479),  we deduce
\begin{align*}
\max_{0\leq i\leq n}I_i\leq& C\, e^{CT}  \E\Big( I_n+\sum_{1\leq i\leq n-1} 
\iti{i-1}{i+1}\normeZ{ Z_r-Z_{t_i}}^2 dr+|\pi|\Big)
\\
\leq
&
 C   \E\Big(  \normeY{\phi\lp X_{T}\rp-\phi\lp X^{\pi}_{T}\rp}^2
+\sum_{1\leq i\leq n}  \iti{i-1}{i}\!\! \lp\normeZ{ Z_r-Z_{t_{i-1}}}^2
+\normeZ{ Z_r-Z_{t_{i}}}^2\rp dr +|\pi|\Big) .
\end{align*}
Since $\phi$ is Lipschitz, Proposition \ref{propositionX} implies that 
$\E|\phi(X_T)-\phi(X_T^\pi)|^2 \leq C |\pi|$; thus  Theorem \ref{L2regularity} implies
\begin{equation}\label{cinqcinq}
  \max_{0\leq i\leq n}\E\normeY{ Y_{t_i}-\schemeZhangY{i}}^2
\leq C|\pi| . 
\end{equation}
Moreover, summing both sides of (\ref{cinqquatre}) over $i$ from 1 to $n-1$
and using Corollary \ref{corZ} we obtain:
\begin{align*} 
\sum_{0\leq i\leq n-2} I_{i}+\frac{\delta}{3}\; \E\int_{t_1}^{T} \! \normeZ{ Z^{\pi}_r-Z_r}^2dr \leq
&
 \sum_{1\leq i < n} \lp 1+ C \Delta t_i \rp I_i+C|\pi|\\
& +C\sum_{1\leq i<n}\E\iti{i-1}{i+1}\!\! \normeZ{ Z_r-Z_{t_i}}^2 dr,\\
\leq
&
C|\pi|+ \sum_{1\leq i\leq n-1} \lp 1+ C \Delta t_i \rp I_i.
\end{align*}
Therefore,
\begin{align*}
 I_{0}+\frac{\delta}{3}\; \E\int_{t_1}^{T} \! \normeZ{ Z^{\pi}_r-Z_r}^2dr 
\leq
&
C|\pi|+  I_{n-1}+C\sum_{1\leq i\leq n-1}   \Delta t_i  I_i
 \end{align*}
Since $\delta<1-\kappa\alpha<3$, using  \eqref{cinqcinq} we deduce
\begin{align}
\frac{\delta}{3}\; \E\int_{0}^{T} \! \normeZ{ Z^{\pi}_r-Z_r}^2dr \leq
&
C|\pi|+  \E\int_{t_{n-1}}^{t_{n}}|Z_r^{\pi}-Z_r|^2dr
+C|\pi|   \E\int_{0}^{T}|Z_r^{\pi}-Z_r|^2dr.\label{numero1}
\end{align}
The equations \eqref{BDSDE} and \eqref{SchemeZhang} imply
\begin{align*}
\iti{n-1}{n}(
Z_r^{\pi}-Z_r)dW_r=&(Y^{\pi}_{t_n}-Y_{t_n})-(Y^{\pi}_{t_{n-1}}-Y_{t_{n-1}})\\
& +\int_{t_{n-1}}^{t_n}\lp f\lp t_n,X^{\pi}_{t_n},Y^{\pi}_{t_n},0\rp-f\lp r,X_r,Y_r, Z_r\rp\rp dr\\
&+\int_{t_{n-1}}^{t_n}\lp g\lp t_n,X^{\pi}_{t_n},Y^{\pi}_{t_n},0\rp- g\lp r,X_r,Y_r,Z_r\rp\rp \dBr.
\end{align*}
The $L^2$-isometry, Schwarz's inequality, \eqref{cinqcinq}, Lemma
\ref{lemmaYYZ}, Propositions \ref{propositionForwardX} and \ref{propositionX} 
\begin{align}
\E\iti{n-1}{n}\lv
Z_r^{\pi}-Z_r\rv^2dr\leq&4\E\lv Y^{\pi}_{t_n}-Y_{t_n}\rv^2+4\E\lv Y^{\pi}_{t_{n-1}}-Y_{t_{n-1}}\rv^2\nonumber\\
& +4|\pi|\E\int_{t_{n-1}}^{t_n}\lv f\lp t_n,X^{\pi}_{t_n},Y^{\pi}_{t_n},0\rp-f\lp r,X_r,Y_r, Z_r\rp\rv^2 dr\nonumber\\
&+4\E\int_{t_{n-1}}^{t_n}\lv g\lp t_n,X^{\pi}_{t_n},Y^{\pi}_{t_n},0\rp- g\lp
r,X_r,Y_r,Z_r\rp\rv^2 dr\nonumber\\
\leq& C|\pi|+C|\pi|\sup_{t_{n-1}\leq r\leq t_n}\E\lp |X_r-X_T|^2+
|X^{\pi}_{t_n}-X_T|^2\rp\nonumber\\
& +C|\pi|\sup_{t_{n-1}\leq r\leq t_n}\E\lp|Y_r-Y_T|^2+ |Y^{\pi}_{t_n}-Y_T|^2+ |Z_r|^2\rp\nonumber\\
\leq& C|\pi|.\label{numero2}
\end{align}
For $|\pi|$ small enough, we have $C|\pi|\leq\delta/6$; thus
\eqref{numero1} and \eqref{numero2}
conclude the proof.
\endproof
\newcommand{\RRX}{\R}
\section{A numerical scheme}
In this section we propose a numerical scheme based on the results of the previous sections. 
First of all, given  $x\in\R^d $, $s<t$ we set:
$$
X_t\lp s,x\rp :=x+\lp t-s\rp b\lp x\rp +\sigma \lp x\rp \lp W_t-W_s\rp.
$$
We clearly have  
$ X^{\pi}_{t_i}=X_{t_i}\lp t_{i-1},X^{\pi}_{t_{i-1}}\rp$ for every  $i=1,\dots,n$. 
Then, given a vector $\lp x_0,\dots,x_i;x_{i+1},\dots ,x_n\rp \in\R^{(i+1)d}\times\R^{n-i}$, set  
$\finX{n+1}=\emptyset$ and for $i=0,\dots,n-1$, let
$$
\vectorX{i}:=\lp x_0,...,x_i\rp , \quad \finX{i+1}:=\lp x_{i+1},...,x_n\rp.
$$
Define by induction, the functions  
 $u^{\pi}_i$, $v^{\pi}_i:\R^{(i+1)d}\times\R^{n-i}\rightarrow\R$ (resp. the random variables
 $U^{\pi}_i,V^{\pi}_i:\R^{(i+1)d}\times\Omega\times\R^{n-i-1}\rightarrow\R$) as follows:
\[ 
 u^{\pi}_n\lp x_0,...,x_n\rp :=\phi\lp x_n\rp ,\quad 
v^{\pi}_n\lp x_0,...,x_n\rp :=0,  
\]
and for $i=0, \dots,  n-1$ let 
\begin{align}
 U^{\pi}_i\lp \vectorX{i},\omega,\finX{i+2}\rp 
:=
&
u^{\pi}_{i+1}\lp \vectorX{i},X_{t_{i+1}}\lp t_{i},x_{i}\rp,\finX{i+2} \rp 
\label{numeroun}\\
&
+f\Big( t_{i+1},X_{t_{i+1}}\lp t_{i},x_{i}\rp ,u^{\pi}_{i+1}\lp \vectorX{i},
X_{t_{i+1}}\lp t_{i},x_{i}\rp,\finX{i+2} \rp ,
\nonumber\\
&
\qquad\qquad v^{\pi}_{i+1}\lp \vectorX{i},X_{t_{i+1}}\lp t_{i},x_{i}\rp,\finX{i+2} \rp  \Big)\Delta t_{i+1},
\nonumber\\
 V^{\pi}_i\lp \vectorX{i},\omega,\finX{i+2}\rp 
:=
&
g\Big( t_{i+1},X_{t_{i+1}}\lp t_{i},x_{i}\rp ,u^{\pi}_{i+1}\lp
 \vectorX{i},X_{t_{i+1}}\lp t_{i},x_{i}\rp,\finX{i+2} \rp ,
\nonumber\\
&
\qquad\qquad v^{\pi}_{i+1}\lp \vectorX{i},X_{t_{i+1}}\lp t_{i},x_{i}\rp,\finX{i+2} \rp  \Big), 
\label{numerodeux}\\
 u^{\pi}_{i}\lp \vectorX{i};\finX{i+1}\rp 
:=
&
\E U^{\pi}_i\lp \vectorX{i},\omega,\finX{i+1}\rp 
+x_{i+1}\E V^{\pi}_i\lp \vectorX{i},\omega,\finX{i+1}\rp , 
\label{numerotrois}\\
v^{\pi}_{i}\lp \vectorX{i};\finX{i+1}\rp 
:=
&
\frac{1}{\Delta t_{i+1}}\E \lp U^{\pi}_i\lp \vectorX{i},\omega,\finX{i+1}\rp \Delta W_{t_{i+1}}\rp \nonumber
\\&\qquad \qquad 
+\frac{x_{i+1}}{\Delta t_{i+1}}
\E \lp V^{\pi}_i\lp \vectorX{i},\omega,\finX{i+1}\rp  \Delta W_{t_{i+1}}\rp. 
\label{numeroquatre}
\end{align}
\begin{theorem}
 We have for all $i=0,...,n$
\begin{align}
Y^{\pi}_{t_i}&=u_i^{\pi}\lp X_{t_0}^{\pi},...,X_{t_i}^{\pi};\Delta B_{t_{i+1}},...,\Delta B_{t_n}\rp \label{Yu}, \\
Z^{\pi,1}_{t_i}& 
=v_i^{\pi}\lp X_{t_0}^{\pi},...,X_{t_i}^{\pi};\Delta B_{t_{i+1}},...,\Delta B_{t_n}\rp . \label{Zv}
\end{align}
\end{theorem}
\proof
We proceed by backward induction. For $i=n$,  by definition $Y^{\pi}_{t_n}=\phi\lp X^{\pi}_{t_n}\rp $, so (\ref{Yu})
 and (\ref{Zv}) hold trivially.\\
Suppose that the result is true for $j=n, n-1, \cdots, i$. The scheme described in \eqref{SchemeZhang}
implies that 
\begin{equation}
Y^{\pi}_{t_{i-1}}=Y^{\pi}_{t_i}+f\lp t_i,\Theta^{\pi,1}_{t_i} \rp
\Delta t_i 
 +g\lp t_i, \Theta^{\pi,1}_{t_i} \rp 
\Delta B_{t_i}-\iti{i-1}{i} Z^{\pi}_rdW_r.\label{Yipi}
\end{equation}
To prove (\ref{Yu}), we take the conditional expectation of (\ref{Yipi}) with respect to 
$\Fchapeau_{t_{i-1}}=\F^W_{\tinitial,t_{i-1}}\vee\F^B_{\tinitial,T}$; this yields
\begin{align*}
 \E\lp  Y^{\pi}_{t_{i-1}}\left|\Fchapeau_{t_{i-1}}\right.\rp =&\E\lp  Y^{\pi}_{t_i}|\Fchapeau_{t_{i-1}}\rp 
+\E\lp  f\lp t_i,\Theta^{\pi,1}_{t_i} \rp\Delta t_i\left|\Fchapeau_{t_{i-1}}\right.\rp \\
& +\E\lp  g\lp t_i,\Theta^{\pi,1}_{t_i} \rp \Delta B_{t_i}\left|\Fchapeau_{t_{i-1}}\right.\rp
 -\E\lp  \iti{i-1}{i} Z^{\pi}_rdW_r\left|\Fchapeau_{t_{i-1}}\right.\rp .
\end{align*}
Using the fact that 
 $\iti{i-1}{i} Z^{\pi}_rdW_r$  is orthogonal to any
 $\Fchapeau_{t_{i-1}}$-measurable random variable, and the induction
 hypothesis we deduce:
\begin{align*}
Y^{\pi}_{t_{i-1}}=&\E\lp  Y^{\pi}_{t_i}\left|\Fchapeau_{t_{i-1}}\right.\rp
 +\E\lp  f\lp t_i,\Theta^{\pi,1}_{t_i} \rp\Delta t_i\left|\Fchapeau_{t_{i-1}}\right.\rp 
 +\E\lp  g\lp t_i, ,\Theta^{\pi,1}_{t_i} 
\rp \left|\Fchapeau_{t_{i-1}}\right.\rp \Delta B_{t_i} .\\
=
&
\E\lp  u_i^{\pi}\lp X_{t_0}^{\pi},...,X_{t_{i-1}}^{\pi},X_{t_i}(t_{i-1},X^{\pi}_{t_{i-1}}),\Delta B_{t_{i+1}},...,
\Delta B_{t_n}\rp  \Big|\Fchapeau_{t_{i-1}}\rp
\\ 
&
+\Delta t_i\;\E\Big(   f\Big( t_i,X_{t_i}(t_{i-1},X^{\pi}_{t_{i-1}}),
u_i^{\pi}\lp  X_{t_0}^{\pi},...,X_{t_{i-1}}^{\pi},X_{t_i}(t_{i-1},X^{\pi}_{t_{i-1}}),\Delta B_{t_{i+1}},...,\Delta B_{t_n}\rp,
\\
&
\qquad\qquad\qquad\qquad v_i^{\pi}\lp X_{t_0}^{\pi},...,X_{t_{i-1}}^{\pi},X_{t_i}(t_{i-1},X^{\pi}_{t_{i-1}}),
\Delta B_{t_{i+1}},...,\Delta B_{t_n}\ \rp \Big)  \Big|\Fchapeau_{t_{i-1}}\Big)  
\\
&
 +\Delta B_i\;\E\Big(  g\Big(  t_i,X_{t_i}(t_{i-1},X^{\pi}_{t_{i-1}}),u_i^{\pi}\Big( X_{t_0}^{\pi},...,
X_{t_{i-1}}^{\pi},X_{t_i}(t_{i-1},X^{\pi}_{t_{i-1}}),\Delta B_{t_{i+1}},...,\Delta B_{t_n}\Big) ,  
\\
& \qquad\qquad\qquad\qquad v_i^{\pi}\lp X_{t_0}^{\pi},...,X_{t_{i-1}}^{\pi},X_{t_i}(t_{i-1},X^{\pi}_{t_{i-1}}),
\Delta B_{t_{i+1}},...,\Delta B_{t_n}\ \rp \Big)\Big|\Fchapeau_{t_{i-1}}\Big).
\end{align*}
Since 
all $\Delta B_{t_l}$, $l=1, \dots, n$ and $X^{\pi}_{t_k}, k=0,\dots,i-1$ are $\widehat{\F}_{t_{i-1}}$
measurable while $W_{t_i}-W_{t_{i-1}}$ is independent of
$\widehat{\F}_{t_{i-1}}$; we deduce (\ref{Yu}).\\
To prove (\ref{Zv}), multiply (\ref{Yipi}) by $\Delta W_{t_i}=W_{t_i}-W_{t_{i-1}}$
and take the conditional expectation with respect to $\Fchapeau_{t_{i-1}}$, this yields 
\begin{align*}
\E\Big(  Y^{\pi}_{t_{i-1}} & \Delta W_{t_i}\Big| \Fchapeau_{t_{i-1}}\Big)    
=
\E\Big( Y^{\pi}_{t_i}\Delta W_{t_i}\Big| \Fchapeau_{t_{i-1}}\Big) 
+\E\lp\left. f\lp t_i,\Theta^{\pi,1}_{t_i} \rp\Delta t_i\Delta W_{t_i}\rv\Fchapeau_{t_{i-1}}\rp  
\\
&
+\E\lp g\lp t_i, \Theta^{\pi,1}_{t_i}
  \rp \Delta B_{t_i}\Delta W_{t_i}\lv\Fchapeau_{t_{i-1}}\right.\rp 
-\E\lp\left. \Delta W_{t_i}\iti{i-1}{i} Z^{\pi}_rdW_r\rv\Fchapeau_{t_{i-1}}\rp .
\end{align*}
Since $Y^{\pi}_{t_{i-1}}\in\Fchapeau_{t_{i-1}}$ and 
$\Delta W_{t_i}$ is independent of $\Fchapeau_{t_{i-1}}$ and centered we deduce 
$$
\E\lp\left.  Y^{\pi}_{t_{i-1}}\Delta W_{t_i}\rv\Fchapeau_{t_{i-1}}\rp   =0 . 
$$
Furthermore, 
\begin{align*}
\E\Big( \Delta W_{t_i} \int_{t_{i-1}}^{t_i} Z^\pi_r dW_r \Big| \widehat{\F}_{t_{i-1}}\Big) =&
\E\Big( \int_{t_{i-1}}^{t_i} Z^\pi_r dr \Big| \widehat{\F}_{t_{i-1}}\Big)\\
=& \E\Big( \int_{t_{i-1}}^{t_i} Z^\pi_r dr \Big| {\F}_{t_{i-1}}\Big)= \Delta t_i Z^{\pi,1}_{t_i}.
\end{align*}
this  completes the proof of (\ref{Zv}).
\endproof
\subsection*{Acknowledgments:} The author wishes  to thank Annie Millet for
helpful comments and  for her precious help in the final preparation of this paper.

\end{document}